\documentclass[11pt]{article}

 \usepackage{amsmath,amssymb,amsthm}
    \usepackage[mathscr]{eucal}
\usepackage{pictexwd,dcpic}

\usepackage{fullpage}

\title{Von Neumann Categories}

\author{\begin{tabular}[t]{c}
        Richard Blute\thanks{Research supported in part by
NSERC.} \,\,\,\,\,\,\,\,\,\,\,\,\,\,Marc Comeau$^*$\\
        {\small Department of Mathematics and Statistics}\\ [-4pt]
        {\small University of Ottawa}\\  [-4pt]
        {\small Ottawa, Ontario, Canada}\\ [-4pt]
        {\small \textsf{
rblute@uottawa.ca\,\,\,\,\,\,\,\,\,\,\,marchy2@hotmail.com}}
\end{tabular}}

\newtheorem{thm}{Theorem}[section]
\newtheorem{proposition}[thm]{Proposition}
\newtheorem{cor}[thm]{Corollary}

\newtheorem{lem}[thm]{Lemma}
\newtheorem{defn}[thm]{Definition}
\newtheorem{example}[thm]{Example}
\newcommand{\rarr}{\rightarrow}
\newcommand{\cC}{\mbox{$\cal C$}}
\newcommand{\cB}{\mbox{$\cal B$}}
\newcommand{\sH}{\mbox{$\sf H$}}
\newcommand{\sK}{\mbox{$\sf K$}}
\newcommand{\ox}{\otimes}

\def\pushright#1{{
    \parfillskip=0pt            
    \widowpenalty=10000         
    \displaywidowpenalty=10000  
    \finalhyphendemerits=0      
    \leavevmode                 
    \unskip                     
    \nobreak                    
    \hfil                       
    \penalty50                  
    \hskip.2em                  
    \null                       
    \hfill                      
    {#1}                        
    \par}}                      

 \def\qed{\pushright{$\Box$}\penalty-700 \smallskip}

\newenvironment{prf}[1]{\begin{trivlist} \item[{\bf ~Proof}#1.]}%
{\qed\end{trivlist}}

\begin{document}

\maketitle

\begin{abstract}
In this paper, we introduce the notion of a {\it von Neumann category}, as a
generalization and categorification of von Neumann algebra. A von Neumann category is a 
premonoidal category with compatible dagger structure which embeds as a double commutant into
a suitable premonoidal category of Hilbert spaces.

The notion was inspired by algebraic quantum field theory. In AQFT, one assigns to
open regions in Minkowski space a $C^*$-algebra, called the local algebra. The local algebras are patched together
to form a global algebra associated to the AQFT. The 
key relativistic assumption is {\it Einstein Causality}, which says that the algebras associated 
to spacelike separated regions commute in the global algebra. 
Premonoidal categories provide a natural framework for lifting such structure from algebras
to categories. Thus von Neumann categories serve as a basis for extending the
abstract quantum mechanics of Abramsky and Coecke to include relativistic effects.

In this paper, we focus on the structure of von Neumann categories. After
giving the basic definitions and examples, we consider constructions typically
associated to von Neumann algebras, and examine their extensions to the category
setting. In particular, we present a crossed product construction for $*$-premonoidal
categories.
\end{abstract}

\section{Introduction}

Algebraic quantum field theory (AQFT) is a mathematically rigorous framework
for modelling the interaction of quantum mechanics in its $C^*$-algebra
interpretation and relativity, as modelled in Minkowski space. It is also
explicitly category-theoretic; essentially an AQFT is a well-behaved functor.
We recommend \cite{Hal,Rob} as references.

One considers Minkowski space
as an ordered set with the causal ordering \cite{Pen}. Then one takes the set of
{\it double cones} or {\it intervals},
that is to say sets of the form:

\[[a,b]=\{x|a\leq x\leq b\}\]
 
Intervals form a partially-ordered set under inclusion. An AQFT is then an
assignment of a $C^*$-algebra to each interval. So we have a map:

\[\cal{U}\mapsto\cal{A}(U)\]

The algebras $\cal{A}(U)$ are called {\it local algebras}. They are the algebras of observables local 
to that region. There are a number of properties, but two that are of interest to us:

\begin{itemize}
\item  The local algebras satisfy that if $\cal{U\subseteq V}$, then $\cal{A(U)\subseteq A(V)}$, i.e. 
the assignment ${\cal A}$ is functorial.
\end{itemize}

For the second condition, note that the set of double cones in Minkowski space is directed, thus one
 can form the directed colimit of the local algebras. The result is
  denoted $\hat{\cal{A}}$, and called the {\it quasilocal algebra}.  The second condition is then:

\begin{itemize}
\item (Einstein Causality) If ${\cal U}$ and ${\cal V}$ are spacelike separated regions, then the local 
algebras 
$\cal{A}(U)$ and $\cal{A}(V)$ pairwise commute in the quasilocal algebra.
\end{itemize}

Einstein causality is the main relativistic assumption, stating that there can be no influence propagated between 
spacelike separated regions. There are typically other axioms, for example involving an action of the
 Poincar\'e group, but these will not concern us here. 

The second influence on this work is the {\it abstract quantum mechanics} of Abramsky and 
Coecke \cite{AC,AC2}. There, quantum mechanics is reformulated away from the notion of $C^*$-algebras and expressed in abstract, categorical terms.  The categorical structure in question is that of 
a  {\it compact closed dagger category}. (See also \cite{Sel} where the structure of these categories is 
examined through the development of a graphical language, and an abstract form of {\it completely 
positive map} is introduced.) The authors of \cite{AC} show that much of the classical theory of 
quantum mechanics can be carried out in this more abstract setting. The authors show for example 
that compact closed dagger categories provide sufficient structure to 
model protocols such as quantum teleportation or entanglement swapping \cite{CN}. The 
correctness of the interpretation basically amounts to the coherence equations of the theory. 
The canonical example of a compact closed dagger category is the category of finite-dimensional Hilbert spaces. Indeed, Selinger has recently shown \cite{Sel2} that the category of finite-dimensional Hilbert spaces is complete for this theory in the sense that an equation follows from the axioms of compact closed dagger
categories if and only if it holds in finite-dimensional Hilbert spaces.
Thus the Abramsky-Coecke axiomatization, while more abstract, is clearly an appropriate level of generality.

But this encoding of teleportation does not take into account that fact that teleportation takes place in 
spacetime. In quantum teleportation, for example, the two participants must pass a classical message. So when
this occurs, they cannot be spacelike 
separated. We believe that an appropriate modification of AQFT would allow for such modelling. More 
specifically,  one should associate some sort of category of local protocols to each region in spacetime.  But what 
structure should the category have? A reasonable first guess would be that of a compact closed 
dagger category. But this leaves open the question of how to express Einstein Causality. 

We propose here modifying the usual notion of compact closed dagger category by replacing the 
monoidal structure with premonoidal structure, as introduced by Power and Robinson \cite{PR}. One 
of the fundamental aspects of monoidal structure in a category is the bifunctoriality of the tensor 
product. That is precisely what is weakened in the definition of premonoidal category. We claim that 
the usual bifunctoriality equation

\[(A\ox f)(g\ox B)=(g\ox B)(A\ox f)\]

\noindent (which is of course then denoted $f\ox g$) can be used to capture the Einstein causality 
condition.

The premonoidal version of AQFT has been developed in the thesis of the second author \cite{Com}. (We note that an alternative approach to modelling this issue has been proposed by Coecke and Lal \cite{CL}. They consider categories in which the tensor product is only partially defined on arrows.)
In this paper, we develop the associated abstract categorical structure. We will view a $*$-premonoidal category as a multiobject version of a $*$-algebra. So we would like to  examine several $C^*$-algebraic constructions. 

Of course, a {\it von Neumann algebra} \cite{Sun} is defined as a subset of $\cB(\sH)$ equal to its own
 double-commutant.  The premonoidal analogue of $\cB(\sH)$ is a category we denote by ${\sf Hilb}
_{\sH}$, where ${\sf H}$ is fixed. Its objects are Hilbert spaces, and an arrow  $f\colon \sK_0\rarr 
\sK_1$ in ${\sf Hilb}_{\sH}$ is an arrow
$f\colon \sK_0\ox \sH\rarr \sK_1\ox \sH$ in ${\sf Hilb}$, the category of (arbitrary) Hilbert spaces and bounded linear maps. The tensor product on objects is the same as 
the tensor product in  ${\sf Hilb}$. The premonoidal structure on arrows is defined below. 

Then, given a set ${\sf A}$ of arrows in  ${\sf Hilb}_{\sH}$, we define its {\it commutant} ${\sf A}'$ to be 
all those arrows of ${\sf Hilb}_{\sH}$ which satisfy the bifunctoriality equation with respect to all arrows 
in ${\sf A}$. We always have that ${\sf A}'$ is a premonoidal subcategory of ${\sf Hilb}_{\sH}$ and ${\sf 
A}\subseteq
{\sf A}''$. Then a {\it von Neumann category} is a subcategory of ${\sf Hilb}_{\sH}$ equal to its own double-commutant. It is always a premonoidal dagger category.

Of interest are constructions for building von Neumann categories. One of the most important ways to build von Neumann algebras is the {\it crossed product} construction \cite{JS}. We present a lifting of this construction to our premonoidal setting. The result is always a  von Neumann category.

\smallskip


\section{Premonoidal categories}

While we assume knowledge of monoidal category theory \cite{CWM}, we here develop the notion of {\it premonoidal category}, due to Power and Robinson \cite{PR}.  As mentioned in the introduction, the intuition behind the definition is that these are monoidal categories, but without assuming the bifunctoriality of the tensor. So in particular, every monoidal category is premonoidal. 
To formalize this intuition, we first introduce {\it binoidal} categories. 

\begin{defn}\emph{
    A {\em binoidal} category consists of a category $\cal{C}$ and functors
$H_{B} : \cal{C}\longrightarrow \cal{C}$ and $K_{B} : \cal{C} \longrightarrow
\cal{C}$ for all objects $B$ in
    $\cal{C}$ and satisfying $H_{B}(C) = K_{C}(B)$ for all pairs of objects $B$.}
\end{defn}

In a binoidal category the object $H_{B}(C)= K_{C}(B)$ is denoted $B\otimes C$ and
for any arrow $f:X\longrightarrow Y$ we write $B\otimes f$ for $H_{B}(f)$ and
$f\otimes B$ for $K_{B}(f)$.  Thus in this new
notation $H_{B} = B\otimes -$ and $K_{B} = - \otimes B$.  Notice that $- \otimes -$
is only a functor when one
of the arguments is fixed, i.e. it is not assumed to be a bifunctor.

\begin{defn}{\em 
    If $\cal{C}$ is a binoidal category and $f : A \longrightarrow C$ 
    and  $g : B \longrightarrow D$ are arrows, define 
 \[g\rtimes f=(g\ox C)(B\ox f)\mbox{\,\,\,\,\,\,\,\,\,\,and \,\,\,\,\,\,\,\,\,\,\,}g\ltimes f=(D\ox f)(g\ox A).\]
 
 Then we say that $f$ is {\em central} if for all arrows $g$, 
 
 \[g\rtimes f=g\ltimes f\mbox{\,\,\,\,\,\,\,\,\,\,and\,\,\,\,\,\,\,\,\,\,}f\rtimes g=f\ltimes g\]

}\end{defn}

\begin{defn}\emph{
    If $\cal{C}$ is a binoidal category and $G, H : \cal{B} \longrightarrow
\cal{C}$ are functors
    then a natural transformation $\alpha : G \Longrightarrow H$ is
{\em central} if its components
    $\alpha_{B}: G(B) \longrightarrow H(B)$ are central maps in $\cal{C}$.}
\end{defn}

\begin{defn}\emph{
A {\em premonoidal category} consists of a binoidal category $\cal{C}$
together with a distinguished
object $I \in |\cal{C}|$ and central natural isomorphisms $\alpha$, $\lambda$
and $\rho$ with components
$\alpha_{A,B,C} : (A\otimes B) \otimes C \longrightarrow A\otimes ( B \otimes C)$,
$\lambda_{A} : I\otimes A
\longrightarrow A$, and $\rho_{A} : A\otimes I \longrightarrow A$.  These
structural isomorphisms must satisfy the same coherence equations as in the
definition of a monoidal category \cite{CWM}.
A premonoidal category is {\em symmetric} if there is a central natural isomorphism
$\tau\colon A\ox B\rarr B\ox A$, again satisfying the usual equations.
}\end{defn}


\subsection{Examples}

\begin{example}\emph{
    If $M$ is a monoid, then $M$ is a one-object premonoidal category. We denote it $M[1]$. It is monoidal if and only if 
    $M$ is abelian. Note that the standard result that in a monoidal category with tensor unit $I$, that the monoid $Hom(I,I)$ is abelian, is false for premonoidal categories. }
\end{example}
\begin{example}\label{E:premFunc}\emph{
    If $\cal{D}$ is any category then define a new category $\cal{C} =
\{{\cal D},{\cal D}\}$
    whose objects are functors $F : \cal{D} \longrightarrow \cal{D}$ and an
arrow $h : F \longrightarrow
    G$ is a {\em transformation,} i.e. consists of arrows $h_{D} : FD
\longrightarrow GD$ for each $D \in
    |\cal{D}|$.  Then $F\otimes G = F\circ G$ for $F,G \in |\cal{C}|$ and
for any transformation $h : F
    \longrightarrow G$ define $(H\otimes h)_{D} = H(h_{D})$ and $(h\otimes H)_{D} =
h_{HD}$.  Then $\cal{C}$
    is a premonoidal category.  If one restricts to transformations which
are natural then one obtains a
    subcategory of $\cal{C}$ which is monoidal.}
\end{example}
\begin{example}\emph{
    Every monoidal category is a premonoidal category.}
\end{example}
\begin{defn}\emph{
    If $\cal{C}$ is a premonoidal category then the {\em centre} of
$\cal{C}$ is the category
    $\cal{Z}(\cal{C})$ with objects the same as those of $\cal{C}$ and
its arrows are the central
    maps in $\cal{C}$.}
\end{defn}

\begin{example}\emph{
    If $M$ is a group, then ${\cal Z}(M[1])$ is just the centre of $G$ viewed as a one-object monoidal
category.}
\end{example}

\begin{proposition}
    The centre $\cal{Z}(\cal{C})$ of a premonoidal category $\cal{C}$
is a monoidal category.
\end{proposition}

\begin{example}{\em 
    Let ${\cal C}$ be a symmetric monoidal category with symmetry $\tau_{X,Y} :
X\otimes Y \longrightarrow
    Y\otimes X$ and let $S\in |{\cal C}|$ be a fixed object.  Define a new
category ${\cal C}_{S}$ as
    follows, the objects are the same as those of $\cal{C}$ and
${\cal C}_{S}(X,Y) = {\cal C}(X\otimes S,
    Y\otimes S)$.  For $Z \in |{\cal C}_{S}|$ and $f \in {\cal C}_{S}(X,Y)$
define $Z\otimes f \in
    {{\cal C}_S}(Z\otimes X, Z\otimes Y)$ as $id_Z\ox f\colon Z\ox X\ox S\rarr Z\ox Y\ox S$ in ${\cal C}$.
    For $f\ox Z$, define it in ${\cal C}$ as: 
    
    \[X\ox Z\ox S\stackrel{\tau_{12}}{\longrightarrow}Z\ox X\ox S\stackrel{id_Z\ox f}{\longrightarrow}
    Z\ox Y\ox S\stackrel{\tau_{12}}{\longrightarrow}Y\ox Z\ox S\]
    
\noindent  The structural isomorphisms for associativity and units come from the corresponding maps in
    ${\cal C}$.
  }
\end{example}

It is straightforward to verify that arrows $g\in {{\cal C}_S}(X,Y)$ of the form 
$g=h\ox id_S$ with $h\colon  X\rarr Y$ in ${\cal C}$ are central. In many examples, all central maps are of this form. For example, we have the following, which was proved in \cite{Com}.

\begin{thm}\emph{
    If $\sH$ is a Hilbert space with $dim (\sH) \geq 1$ then
${\cal Z}(\mathbf{Hilb}_{\sH}) \simeq \mathbf{Hilb}$. More precisely, all central maps are of the form $f=\hat{f}\ox id_H\colon
K_1\ox H\rarr K_2\ox H$, with $\hat{f}\colon K_1\rarr K_2$. 
}\end{thm}

\begin{prf}{}
    If $dim(\sH) = 1 $ then $ \sH \cong \mathbb{C}$ and clearly $\mathbf{Hilb}_{\sH} \simeq \mathbf{Hilb}$.  So now
    suppose that $dim(\sH) > 1$ and that $f \in \mathbf{Hilb}_{\sH}(X,Y)$ is central.  Then we will show that $f =
    \hat{f}\otimes id_{\sH}$ for some bounded linear map $\hat{f}: X \longrightarrow Y$.  Let ${\cal B}_{\sH} = \{ h_{j}
    \mid j \in J \}$ be an orthonormal basis for $\sH$.  Then for $a \neq b \in J$ define $T_{a,b} : \sH\otimes \sH
    \longrightarrow \sH \otimes \sH$ by $T_{a,b}(h_{i}\otimes h_{j} ) = (\delta_{i,a}\delta_{j,b} +
    \delta_{i,b}\delta_{j,a}) h_{j}\otimes h_{i}$ where $\delta_{p,q} = 1$ if $p = q$ and $\delta_{p,q} = 0$
    otherwise.  
    
    Now notice that the vector subspace $(\sH\otimes \sH)_{a,b} = \{ \lambda h_{a}\otimes h_{b} +
    \mu h_{b}\otimes h_{a} \mid \lambda, \mu \in \mathbb{C} \}$ is a finite-dimensional linear subspace of
    $\sH\otimes \sH$ and hence is closed.  Moreover the map $T_{a,b}$ is  then just the projection onto the closed
    subspace $(\sH\otimes \sH)_{a,b}$ followed by a twist and is therefore continuous.

    Now suppose that ${\cal B}_{X} = \{ e_{i} \mid i\in I
    \}$ and ${\cal B}_{Y} = \{ g_{k} \mid k\in K \}$ are orthonormal bases for $X$ and $Y$ respectively.  We now compute
    $f\rtimes T_{a,b} : X \otimes \sH\otimes \sH \longrightarrow Y\otimes \sH\otimes \sH$ on basis elements.
    \begin{align}
        f\rtimes T_{a,b} ( e_{i}\otimes h_{j}\otimes h_{k}) &=(f\otimes \sH)(X\otimes T_{a,b} (e_{i}\otimes h_{j}\otimes
                                                            h_{k}))\\ \notag
                                                            &=(f\otimes \sH)[(\delta_{j,a}\delta_{k,b} +
                                                            \delta_{j,b}\delta_{k,a}) e_{i}\otimes h_{k} \otimes
                                                            h_{j}]
    \end{align}
    Note that when restricting to the diagonal $j = k = a$ and $a\neq b$, we get that
    
\[ f\rtimes T_{a,b}(e_{i}\otimes h_{j}\otimes h_{k})= f\rtimes T_{a,b}(e_{i}\otimes h_{a}\otimes h_{a})
= 0.\]
    
    On the other hand we now calculate $f\ltimes T_{a,b}$ applied to $( e_{i}\otimes h_{j}\otimes h_{k})$.
    First observe that
    \begin{align}
        f\otimes \sH(e_{i}\otimes h_{j}\otimes h_{k}) &= (\tau\otimes id_{\sH})(h_{j}\otimes f(e_{i} \otimes
                                                    h_{k})) \\ \notag
                                                    &= (\tau\otimes id_{\sH}) h_{j}\otimes (\sum_{r\in K, \,
                                                    p \in J} c_{i,k}^{r,p}g_{r}\otimes h_{p})\\ \notag
                                                    &= \sum_{r \in K,\, p \in J} c_{i,k}^{r,p}g_{r}\otimes h_{j}
                                                    \otimes h_{p}.
    \end{align}
    Therefore:
    \begin{align}
        f\ltimes T_{a,b}( e_{i}\otimes h_{j}\otimes h_{k}) &= (Y\otimes T_{a,b})(f\otimes \sH)( e_{i}\otimes
                                                            h_{j}\otimes h_{k})\\ \notag
                                                            &= \sum_{r \in K, \, p \in J}c_{i,k}^{r,p}
                                                            (\delta_{j,a}\delta_{p,b} +
                                                            \delta_{j,b}\delta_{p,a}) g_{r} \otimes h_{p}
                                                            \otimes h_{j}.
    \end{align}
 Restricting to the same diagonal with $j = k = a$, the above becomes
    \begin{align}
        \sum_{r\in K, \, p \in J}c_{i,a}^{r,p}\delta_{p,b}g_{r}\otimes h_{p}\otimes h_{a} = \sum_{r \in
        K}c_{i,a}^{r,b}g_{r} \otimes h_{b}\otimes h_{a}.
    \end{align}
   As $f$ is central it follows that
    \[
        \sum_{r \in K}c_{i,a}^{r,b}g_{r} \otimes h_{b}\otimes h_{a} = 0.
    \]
    So we have shown that $ c_{i,a}^{r,b} = \delta_{a,b} c_{i,a}^{r,b}$,
 for all $r \in K$.  By a similar calculation, one shows that $c_{i,a}^{r,a} = c_{i,b}^{r,b}$ for 
 all $a, \, b \in J$.

    Now fix $a \in J$ and define $d_{i}^{r} = c_{i,a}^{r,a}$ for all $r \in K$.  We 
    have
   
\[        f(e_{i}\otimes h_{a}) =\sum_{r \in K} d_{i}^{r} g_{r} \otimes h_{a}
                                =(\sum_{r \in K} d_{i}^{r} g_{r})\otimes h_{a}
                                =\widehat{f}(e_{i})\otimes h_{a}.\]

    The map $\widehat{f}$ is defined by the equation $\widehat{f}(e_{i})= \sum_{r \in K} d_{i}^{r}g_{r}$.
    Note that $\widehat{f}$ is bounded, since $f=\widehat{f}\ox id$ is bounded.
    
  \end{prf}

 \section{von Neumann categories}

\subsection{Premonoidal $*$-structure}

Much of the theory of $*$-structures on categories appears in \cite{GLR} as a part of their program of analyzing Haag's definition of algebraic quantum field theory \cite{Hal}. It was a crucial component of the Doplicher-Roberts reconstruction theorem \cite{DR1}. 

In this section, we extend this theory to the premonoidal setting. It is straightforward to see that for a monoidal category, viewed as premonoidal, the definitions are the same. 

\begin{defn}\emph{
\begin{itemize}
 \item   An {\em Ab-premonoidal category} is a premonoidal category $\cal{C}$
such that for all objects $X,Y$
    in $\cal{C}$ the set ${\cal C}(X,Y)$ is equipped with an abelian group
structure.  Moveover if $f, g
    \in {\cal C}(X,Y)$ and $h\in {\cal C}(A,X)$ and $k \in {\cal C}(Y,B)$
then $(f+g)\circ h = f\circ h
    + g\circ h$ and $k\circ (f+g) = k\circ f + k \circ g$.  In addition we also
require that for all objects $A$
    the functions $A\otimes - : {\cal C}(X,Y) \longrightarrow
{\cal C}(A\otimes X, A\otimes Y)$ and $-
    \otimes A : {\cal C}(X,Y) \longrightarrow {\cal C}(X\otimes A, Y \otimes
A)$ are group homomorphisms
    for all objects $X$, $Y$ in ${\cal C}$.
\item   A {\em $\mathbb{C}$-linear premonoidal category} is a a premonoidal category
$\cal{C}$ in which
    every hom-set ${\cal C}(X,Y)$ is a complex vector space and the composition
map $(f,g) \mapsto g\circ f$
    is bilinear and the functions $A\otimes - :{ \cal C}(X,Y) \longrightarrow
{\cal C}(A\otimes X, A\otimes
    Y)$ and $- \otimes A : {\cal C}(X,Y) \longrightarrow {\cal C}(X\otimes A,
Y \otimes A)$ are
    $\sf{C}$-linear for all objects $A$, $X$, and $Y$ in $\cal{C}$.
\item
    A {\em positive $\ast$-operation} on a $\mathbb{C}$-linear premonoidal category
$\cal{C}$ is a family
    of functions assigning to each arrow $s \in {\cal C}(X,Y)$ an arrow
$s^{\ast} \in {\cal C}(Y,X)$ with
    $(g\circ f)^{\ast} = f^{\ast} \circ g^{\ast}$ for composable arrows $f$ and $g$
and $id_{A}^{\ast} = id_{A}$.
    The map $s \mapsto s^{\ast}$ must be antilinear, and satisfy $(s^{\ast})^{\ast}
= s$ and if $s^{\ast} \circ s
     = 0$ then $s = 0$.  \\
    \mbox{\,\,} \\
     We also require 
     that for all arrows $g$ in $\cal{C}$ and objects $A$, that $(A\otimes f)^{\ast} = A\otimes f^{\ast}$ and
    $(f\otimes A)^{\ast} =  f^{\ast} \otimes A$.  Finally we will insist
    that $\alpha ^{\ast} = \alpha^{-1}$, $\lambda ^{\ast} = \lambda ^{-1}$, and
$\rho ^{\ast} = \rho ^{-1}$ and
    in the case of symmetry that $\tau^{\ast} = \tau^{-1}$.  A $\mathbb{C}$-linear
premonoidal category equipped with a positive
    $\ast$-operation is called a {\em premonoidal $\ast$-category}.
 \item  A {\em premonoidal $C^{\ast}$-category} is a premonoidal $\ast$-category $\mathcal{C}$ such that
    $\mathcal{C}(X,Y)$ is a Banach space with norm denoted by $\| \cdot \|$ such that $\|s \circ t
    \| \leq \|s\| \|t\|$ and $\|s^{\ast} \circ s \| = \|s \|^{2}$ and $\| A\otimes s
    \|= \| s\|= \|s \otimes A \|$ for all $ s :
    X \longrightarrow Y$ and $t : Y \longrightarrow Z$ and objects $A$. If the category happens to be monoidal, we say that the category is a $C^{\ast}$-{\em tensor category}.
  \end{itemize}}
\end{defn}

 Again by virtue of the definition of premonoidal $C^{\ast}$-categories we have the following result.
\begin{lem}
If $\mathcal{C}$ is a premonoidal $C^{\ast}$-category then $\mathcal{Z}(\mathcal{C})$ is a $C^{\ast}$-tensor
category.
\end{lem}

\subsection{Commutants in premonoidal categories}

The results of this section are inspired by the theory of {\it von Neumann algebras}. They will be the basis of our definition of {\it von Neumann category} below.


\begin{defn}{\em 
    Let $\mathcal{A}$ be a set of objects and arrows in a $*$-premonoidal category $\mathcal{C}$.  Then the
    {\em commutant} of $\mathcal{A}$, denoted  $\mathcal{A}'$, will be the category with objects the same
    as those of $\mathcal{C}$ and its arrows will be arrows $f : A \longrightarrow B$ in $\mathcal{C}$ such that
    $f \ltimes g = f \rtimes g$ and $g \ltimes f = g\rtimes f$ for all arrows $g$ in $\mathcal{A}$.  }
\end{defn}

\begin{thm}
$\mathcal{A}'$ is a $*$-premonoidal category.
\end{thm}
\begin{prf}{}
    We start by showing that $\mathcal{A}'$ is a category.  For each object $A$ the identity map $id_{A}$ is a
    central map in $\mathcal{C}$ and thus a map in $\mathcal{A}'$, so $\mathcal{A}'$ contains identities.  Next
    we must show that given $f: A \longrightarrow B$, and $ e: B \longrightarrow C$ in $\mathcal{A}'$, then
    $e\circ f : A \longrightarrow C$ is an arrow in $\mathcal{A}'$.  Indeed let $g : X \longrightarrow Y$ be an
    arrow in $\mathcal{A}$.

     \begin{align}\notag
        (e \circ f) \ltimes g &= [C\otimes g][(e\circ f)\otimes X]\\ \notag
                           &= ([C\otimes g][e\otimes X])[f\otimes X]\\ \notag
                           &= [e\otimes Y]([B \otimes g] [f\otimes X])\\ \notag
                         &= ([e\otimes Y][f\otimes Y]) [A \otimes g]\\ \notag
                         &= [e\circ f \otimes Y] [A\otimes g]
                        = (e\circ f) \rtimes g
   \end{align}

    Similarly we can show that $g \ltimes (e \circ f) = g \rtimes (e \circ f)$.  Hence $e \circ f$ is an arrow
    in $\mathcal{A}'$.  Clearly this composition is associative and unital thus $\mathcal{A}'$ is a category.
   
    We now establish the premonoidal structure on the commutant category.  Given objects $A$ and $B$ of
    $\mathcal{A}'$, we define $A\otimes_{\mathcal{A}'} B=A \otimes B$.  If also $f : X
    \longrightarrow Y$ in $\mathcal{A}'$ then we define $A\otimes_{\mathcal{A}'} f = A\otimes f$ and $f
    \otimes_{\mathcal{A}'} A = f \otimes A$.  It is an exercise in diagram chasing that 
    $(A\otimes f) \ltimes g = (A\otimes f) \rtimes g$
    and similarly one can check that $g \ltimes (A\otimes f) = g \rtimes (A\otimes f)$.  
    Hence $(A\otimes f)$ is
    an arrow in $\mathcal{A}'$ and likewise so is $(f\otimes A)$.  Now since $\mathcal{Z}(\mathcal{C}) \subseteq
    \mathcal{A}'$ it follows that the remaining requirements for $\mathcal{A}'$ to be a premonoidal category are
    all satisfied since all the relevant diagrams that must commute are diagrams which live in the centre and
    commute there. 
    
    Finally one must check that if $f\in \mathcal{A}'$, then so is $f^*$, but this follows from the functoriality of the $(-)^*$ and the fact that it commutes with all the relevant structure. 
\end{prf}

\subsection{The definition}

The following definition is our attempt at generalizing the notion of a von Neumann algebra, much like our
definition of premonoidal \(C^{\ast}\)-category is an attempt at generalizing the notion of a
\(C^{\ast}\)-algebra.

\begin{defn}\emph{
Let \(\mathcal{A} \subseteq \mathcal{C}\) be a premonoidal \(C^{\ast}\)-subcategory of a premonoidal
\(C^{\ast}\)-category \(\mathcal{C}\).  Then \(\mathcal{A}\) is called a {\em \(\mathcal{C}\)-von Neumann
category} just in case \(\mathcal{A}''(X,Y) = \mathcal{A}(X,Y)\) for all objects \(X\) and \(Y\) in
\(\mathcal{A}\). When \(\mathcal{C} = \mathbf{Hilb}_{\sH}\) then \(\mathcal{A}\) is simply called a {\em von
Neumann category}.}
\end{defn}

One can imagine looking at a more general notion in which there is no normed structure for example, but we choose to stay with this level of generality for the moment. 
A natural question to ask is whether a one-object von Neumann category is a von Neumann algebra.

\begin{thm}
If \(\mathcal{A} \subseteq \mathbf{Hilb}_{\sH}\) is a von Neumann category, then
\(\mathcal{A}(\mathbb{C},\mathbb{C}) \) has the structure of a von Neumann algebra.
\end{thm}

\begin{prf}{}
Let \(\mathcal{M} =  \mathcal{A}(\mathbb{C},\mathbb{C}) \).  First notice that  \(\mathcal{M}\) is a \(\ast\)-subalgebra of \(\mathfrak{B}(\mathbb{C}\otimes \sH)\cong\mathfrak{B}(\sH)\). 

We will show that \(S \in \mathcal{A}'(\mathbb{C},\mathbb{C})\) if and only if \(S\circ T = T\circ S\) for all \(T \in
\mathcal{M}\).  Indeed \(S\ltimes T = S\rtimes T\) means that the following diagram commutes:

\begin{equation}
    \begindc{\commdiag}[30]
        \obj(1,1)[]{$\mathbb{C}\otimes \mathbb{C} \otimes \sH$}
        \obj(5,1)[]{$\mathbb{C}\otimes \mathbb{C} \otimes \sH$}
        \obj(9,1)[]{$\mathbb{C}\otimes \mathbb{C} \otimes \sH$}
        \obj(13,1)[]{$\mathbb{C}\otimes \mathbb{C} \otimes \sH$}
        \obj(1,3)[]{$\mathbb{C}\otimes \mathbb{C} \otimes \sH$}
        \obj(7,3)[]{$\mathbb{C}\otimes \mathbb{C} \otimes \sH$}
        \obj(13,3)[]{$\mathbb{C}\otimes \mathbb{C} \otimes \sH$}
        \mor(1,3)(7,3)[30,30]{$id_{\mathbb{C}}\otimes S$}
        \mor(7,3)(13,3)[30,30]{$\tau_{\mathbb{C},\mathbb{C}}\otimes id_{\sH}$}
        \mor(13,3)(13,1){$id_{\mathbb{C}}\otimes T$}
        \mor(1,3)(1,1){$\tau_{\mathbb{C},\mathbb{C}}\otimes id_{\sH}$}[\atright, \solidarrow]
        \mor(1,1)(5,1)[30,30]{$id_{\mathbb{C}}\otimes T$}[\atright, \solidarrow]
        \mor(5,1)(9,1)[30,30]{$\tau_{\mathbb{C},\mathbb{C}}\otimes id_{\sH}$}[\atright, \solidarrow]
        \mor(9,1)(13,1)[30,30]{$id_{\mathbb{C}}\otimes S$}[\atright, \solidarrow]
    \enddc
\end{equation}

Now recall that \(\tau_{\mathbb{C},\mathbb{C}} = id\).  Using this fact in the above diagram one gets
\begin{equation}
    (id_{\mathbb{C}}\otimes T) \circ (id_{\mathbb{C}}\otimes S) = (id_{\mathbb{C}}\otimes S) \circ (id_{\mathbb{C}}\otimes
    T)
\end{equation}
and thus \(id_{\mathbb{C}}\otimes (T\circ S) = id_{\mathbb{C}}\otimes (S\circ T)\) and this occurs if and only
if \(T\circ S = S\circ T\).  

We denote the commutant of the {\it algebra}
\(\mathcal{M}\) in \(\mathfrak{B}(\mathbb{C}\otimes H)\) by \(\mathcal{M}'\).  Let \(\mathcal{N} =
\mathcal{A}'(\mathbb{C},\mathbb{C})\).  Note that \(\mathcal{A}\) is a von Neumann category, and hence \(
\mathcal{A}(\mathbb{C},\mathbb{C}) = \mathcal{A}''(\mathbb{C},\mathbb{C})\).  Thus \(\mathcal{M} =
\mathcal{N}'\), and clearly as \(\mathcal{M}' = \mathcal{N}\) it follows \(\mathcal{M}'' = \mathcal{N}' =
\mathcal{M} \) showing that \(\mathcal{M}\) is a von Neumann algebra.
\end{prf}

\begin{cor}\label{Cor:OneObVnc}
Every one-object von Neumann category is a von Neumann algebra.
\end{cor}

Thus the above corollary justifies our claim that a von Neumann category is an appropriate generalization of the
notion of a von Neumann algebra.  Before providing some concrete examples of von Neumann categories we will
first establish some analogues of classical results found in the theory of von Neumann algebras.

\begin{proposition}\label{P:normClosedVnc}
If \(\mathcal{A}\) is a set of objects and arrows in a premonoidal \(C^{\ast}\)-category \(\mathcal{C}\) closed
under \(\ast\), then \(\mathcal{A}'\) is a premonoidal \(C^{\ast}\)-category.  In particular, it is a
\(\mathcal{C}\)-von Neumann category.
\end{proposition}
\begin{prf}{}
We already have that \(\mathcal{A'}\) is a premonoidal \(\ast\)-category. Furthermore
each hom-set \(\mathcal{A}'(X,Y)\) is a normed linear subspace of \(\mathcal{C}(X,Y)\) with norm 
coming from  the \(C^{\ast}\)-structure on \(\mathcal{C}\).  Thus it remains to show that each space \(\mathcal{A}'(X,Y)\) is
complete with respect to its norm. 

Notice that for any arrow \(f :A \longrightarrow B\) the linear map \( \zeta_{f} : \mathcal{C}(C,D)
\longrightarrow \mathcal{C}(A\otimes C,B\otimes D)\) given by \(\zeta_{f}(g)=f\ltimes g - f\rtimes g = (B\otimes
g)\circ (f\otimes C) - (f\otimes D)\circ (A\otimes g)\) is bounded.  Similarly the linear map \(\eta_{f} :
\mathcal{C}(B,D) \longrightarrow \mathcal{C}(C\otimes A,D\otimes B)\) given by 
\(\eta_{f}(g) = g\ltimes f -
g\rtimes f = (D\otimes f) \circ (g \otimes A) - (g\otimes B)\circ (C\otimes f)\) is bounded.  So let \((g_{j})\)
be a cauchy sequence in \(\mathcal{A}'(B,D)\), then by completeness of \(\mathcal{C}(B,D)\) it converges to a map
\(g = \lim g_{j}\) in \(\mathcal{C}(B,D)\).  Now for any arrow \(f: A \longrightarrow C\) in \(\mathcal{A}\) we
have that:

\[\zeta_{f}(g) = \zeta_{f}(\lim g_{j}) 
= \lim \zeta_{f}(g_{j})= 0.\]

Similarly we also have that \(\eta_{f}(g) =0\) for any arrow \(f\) in \(\mathcal{A}\) and thus \(g \in
\mathcal{A}'(B,D)\).  Hence we have shown that \(\mathcal{A}'(B,D)\) is closed, establishing that
\(\mathcal{A}'\) is a premonoidal \(C^{\ast}\)-category.

To see that \(\mathcal{A}'\) is a \(\mathcal{C}\)-von Neumann category we observe that \(\mathcal{A} \subseteq
\mathcal{A}''\) and taking commutants we get \(\mathcal{A}''' \subseteq \mathcal{A}'\).  On the other hand we
also have that \(\mathcal{A}' \subseteq \mathcal{A}'''\) and thus the result follows that \(\mathcal{A}''' =
\mathcal{A}'\).
\end{prf}

\subsection{Examples of von Neumann categories}\label{S:exampsVNC}

At this point we feel that some examples of von Neumann categories are in order.  We will also use this
opportunity to draw further parallels between our theory and the classical one.

\begin{example}\emph{
By Corollary \ref{Cor:OneObVnc}, every von Neumann algebra \(\mathcal{M}\) can be viewed as a one-object von
Neumann category.}
\end{example}
\begin{example}\emph{
If \(\mathcal{C}\) is a premonoidal \(C^{\ast}\)-category, then \(\mathcal{C}\) and \(\mathcal{Z}(\mathcal{C})\)
are \(\mathcal{C}\)-von Neumann categories.  This is clear since \(\mathcal{C} = \mathcal{Z}(\mathcal{C})'\). In
the case \(\mathcal{C} = \mathbf{Hilb}_{\sH}\), we see that \(\mathcal{Z}(\mathbf{Hilb}_{\sH}) \simeq \mathbf{Hilb}\)
is a von Neumann category.  In the case that \(\mathcal{C}\) is a von Neumann algebra viewed as a one-object von
Neumann category, we get that centre of a von Neumann algebra is again a von Neumann algebra.}
\end{example}

The above example motivates the following comparison.  If \(\sH\) is a Hilbert space then
\(\mathfrak{B}(\sH)\) is a von Neumann algebra and the centre of \(\mathfrak{B}(\sH)\) is
 \(\mathbb{C}\).  Now by
the above example, \(\mathfrak{B}(\sH)\) can be viewed as a one-object von Neumann category on
\(\mathbf{Hilb}_{\sH}\) and its centre will be the subcategory with object \(\mathbb{C}\) and
 will have as arrows
the central maps on this object.  Thus we think of \(\mathbf{Hilb}_{\sH}\) as a multi-object version of the
classical \(\mathfrak{B}(\sH)\) and likewise since \(\mathcal{Z}(\mathbf{Hilb}_{\sH}) 
\simeq \mathbf{Hilb}\) we
think of \(\mathbf{Hilb}\) as playing the role of the complex numbers \(\mathbb{C}\).

Continuing on with more examples of von Neumann categories, we will 
consider premonoidal \(C^{\ast}\)-categories
that arise as functor categories.

\begin{example}\label{Ex:maxAbelianVnc}\emph{
Suppose that \(\mathcal{D}\) is a \(C^{\ast}\)-category, let \(\{\mathcal{D},\mathcal{D}\}_{\ast}\) be the
premonoidal category whose objects are \(\ast\)-functors and an 
arrow \(t: F \longrightarrow G\) consists of a
family of maps \(t_{A} : FA \longrightarrow GA \) in \(\mathcal{D}\) such that the set \( \{\|t_{A}\| \}\) is
bounded, call these arrows \emph{bounded transformations}.  We should note that this 
example is a premonoidal
variation on an example of \cite{GLR},  of a \(C^{\ast}\)-category.  Now given a map \(t: F
\longrightarrow G\) then one defines \( \|t\| \equiv \sup_{A}\|t_{A}\| \), which yields a norm on the linear
space \(\{\mathcal{D},\mathcal{D}\}_{\ast}(F,G)\) where addition and scalar multiplication are defined
point-wise.}

\emph{The premonoidal structure on this category is the same as the one described in 
Example \ref{E:premFunc}.
Namely given two \(\ast\)-functors \(F\) and \(G\) we define \(F\otimes G \equiv F\circ G\)  which 
is clearly
again a \(\ast\)-functor.  Further given a transformation \(t: F \longrightarrow G\) and a \(\ast\)-functor
\(H\) then we define \((H\otimes t)_{A} \equiv H(t_{A})\) and \((t\otimes H)_{A} \equiv t_{HA} \).  Now it is
clear that \(\{\|(t\otimes H)_{A}\|\}\) is bounded and since \(\|H(f)\| \leq \|f\|\) for all arrows \(f\) it
follows that \(\{\|(H\otimes t)_{A}\|\}\) is also bounded. Now let \(\mathcal{D}^{\mathcal{D}}\) denote
 the wide
subcategory of \(\{\mathcal{D},\mathcal{D}\}_{\ast}\) whose arrows are the bounded natural 
transformations.
Observing that a constant functor is a \(\ast\)-functor one can show that 
\((\mathcal{D}^{\mathcal{D}})' = \mathcal{D}^{\mathcal{D}}\).  Hence \(\mathcal{D}^{\mathcal{D}}\)
is a \(\{\mathcal{D},\mathcal{D}\}_{\ast}\)-von Neumann category.}
\end{example}

\section{Premonoidal crossed products}

\subsection{Crossed products of von Neumann algebras}

We first review the traditional construction. For more information, see \cite{JS}.

Suppose we have a von Neumann algebra $M$ with a
presentation $M\subseteq\cB(\sH)$. Suppose a discrete group $G$
acts on $M$. The crossed product is a von Neumann
algebra $\tilde{M}$ and embeddings

\[\pi\colon M\rarr\tilde{M}\,\,\,\,\,\,\,\,\,\,\,\,\,\lambda\colon
G\rarr\tilde{M}\]

\noindent such that the image of $G$ consists of unitaries, and the images
are related by the conjugation equation

\[\pi(g\cdot a)=\lambda(g)\pi(a)\lambda(g)^*\]

One begins by building a Hilbert space $\tilde{\sH}$, the square-summable functions from $G$ to
$\sH$:

\[\tilde{\sH}=\{\zeta\colon G\rarr\sH|\sum_{g\in G}||\zeta(g)||^2<\infty\}\]

Then define embeddings into $\cB(\tilde{\sH})$ by

\[[\pi(a)(\zeta)](g)=(g^{-1}\cdot a)(\zeta (g))  \mbox{  \,\,\,\,\,\,\,\,  }[\lambda(g)(\zeta)](u)=\zeta(g^{-1}u)\]

\smallskip

\noindent Then the {\it crossed product von Neumann algebra} is defined as $\tilde{M}\subseteq 
\cB(\tilde{\sH})$.

\[\tilde{M}=[\pi(M)\cup\lambda(G)]''\]

\subsection{Premonoidal version}

Remember that our general program is to categorify by replacing $\cB(\sH)$
with ${\sf Hilb}_{\sH}$, and von Neumann algebras with von Neumann categories.

First, we need the analog of a discrete group action on a von Neumann algebra. So
let $G$ be a discrete group, viewed as a premonoidal category $G[1]$. Let $\cC\subseteq
{\sf Hilb}_{\sH}$ be a von Neumann category.  

There are several reasonable levels of generality
we could consider. For the present paper, we choose the functorial action of $G$ to act as the identity 
on objects of $\cC$. So, given an arrow $f\colon \sK\rarr \sK'$ in $\cC$ and $g\in G$, we have an 
action $g\bullet f\colon \sK\rarr \sK'$, satisfying the evident equations.  

We note that the action of $G$ on $\cC$ in particular induces an action 
$\bullet\colon G\times Hom_{\cC}(I,I)\rarr Hom_{\cC}(I,I)$ and recall that by a previous result, we know that $M=Hom_{\cC}(I,I)$
is a von Neumann algebra. This makes the connection to the traditional crossed product even more 
evident. 

To construct a crossed product, we will use the same Hilbert space $\tilde{\sH}$  defined above. 
But it is easier to use an isomorphic description \cite{JS}. So note that $\tilde{\sH}\cong\sH\ox\ell^2(G)$, where 

\[\ell^2(G)=\{f\colon G\rarr\mathbb{C}\mbox{   such that   }\sum_{g\in G}|f(g)|^2<\infty\}\]

Then an (orthonormal) basis for $\tilde{\sH}$ is given by $\{e_i\ox\delta_g\}$, where the $e_i$'s range over a basis for 
\sH, and the $\delta_g$'s are the basis for $\ell^2(G)$ consisting of the Kronecker deltas.

To construct our crossed product, we will embed  $\cC$ and $G$ into ${\sf Hilb}_{\tilde{\sH}}$, and 
then take the double commutant, as above.

Now we need two operations. First I need $\lambda\colon G[1]\rarr{\sf Hilb}_{\tilde{\sH}}$. Since $G[1]$ is a one-object premonoidal category, and a (strong) premonoidal functor must take the tensor unit to itself, this reduces to a map 
$\lambda\colon G\rarr \cB(\tilde{\sH})$, so we can in fact use the definition from the previous section. With the new basis, this is written as

\[\lambda(g)(e_i\ox \delta_{g'})=e_i\ox \delta_{gg'}\]

Second, we need a premonoidal functor $\pi\colon\cC\rarr {\sf Hilb}_{\tilde{\sH}}$. The functor will be the identity on objects. Suppose we have $f\colon X\rarr Y$ in \cC. Thus  $f\colon X\ox\sH\rarr Y\ox\sH$ in {\sf Hilb}. Then $\pi(f)$ should be an arrow $X\rarr Y$ in ${\sf Hilb}_{\tilde{\sH}}$. Define it by 

\[\pi(f)(x\ox e_i\ox\delta_g)=(g^{-1}\bullet f)(x\ox e_i)\ox\delta_g\]

\begin{defn}{\em
The {\em crossed product} of $G$ and $\cC$ will be the double commutant in ${\sf Hilb}_{\tilde{\sH}}$
of the images of $G[1]$ and $\cC$ under the functors $\lambda$ and $\pi$.
}\end{defn}

We then note that we have the following analog of the usual crossed product equation:

\begin{lem} With notation as above, we have:

\[\pi(g\bullet f)=(id\ox\lambda(g))\circ\pi(f)\circ (id\ox\lambda(g)^*)\]

\end{lem}

\section{Conclusion}

We believe this new notion of von Neumann category opens up the possibility of analysis of a 
categorified theory of von Neumann algebras. At this point it would be quite reasonable to introduce a 
theory of {\it factors} \cite{JS} for von Neumann categories and attempt to classify these. We think the most 
likely definition is that a factor is a von Neumann category whose center is {\sf Hilb}.  Then of course we would need to develop integrals over von Neumann categories, and determine which von Neumann categories so arise.
Also this theory allows us to talk about double commutants in categories not based on Hilbert spaces 
at all. We have already given some examples. It remains to be seen how interesting this theory 
is at that level of generality.

Along these lines, the final section of this paper develops only discrete crossed products. But of course the theory of crossed products of von Neumann algebras is much more general than this \cite{Will}. In particular, one would need to develop a theory of integrals to enrich this idea further. 
On the other hand, we could consider crossed products of more general 
premonoidal categories acting on a von Neumann category.

One of the most significant results in the field of algebraic quantum field theory is the {\it Doplicher-
Roberts reconstruction theorem} which demonstrates that every compact closed $C^*$-category is 
equivalent to the category of finite-dimensional representations of an essentially unique compact 
group. See \cite{DR1} for the original result and \cite{Hal} including its appendix by M\"uger for an alternate proof 
and discussion of the result's significance. We are very interested in the analog of this result in the 
premonoidal setting. This is ongoing work.

An important issue in algebraic quantum field theory is the modelling of {\it open systems}, i.e. those 
which can interact with their environment. The difficulties of this are explored in \cite{CH}. We believe 
that using the category ${\sf Hilb}_{\sH}$ can lead to new insight into this issue. In particular, one may 
choose the Hilbert space \sH\ so as to model the environment. So, for example, a protocol of type 
$X\rarr Y$ would be modelled as a morphism of the form $X\ox\sH\rarr Y\ox\sH$, and so may encode
information about such interaction.


\begin{thebibliography}{9}

\bibitem{AC} S. Abramsky, B. Coecke.
\newblock  A categorical semantics of quantum protocols. {\it in}
Proceedings of the 19th Annual IEEE Symposium on Logic in Computer Science
2004, IEEE Comput. Soc. Press, (2004).

\bibitem{AC2} S. Abramsky, B. Coecke. Physics from computer science. International Journal of 
Unconventional Computing 3, pp. 179-197, (2007).

\bibitem{CH} R. Clifton, H. Halvorson. Entanglement and open systems in algebraic quantum field 
theory,  Philosophy of Science 69, pp.1Ð28, (2002).

\bibitem{CN} I. Chuang, M. Nielsen, {\it Quantum Computation and Quantum
Information}, Cambridge University Press, (2000).

\bibitem{CL} B. Coecke, R. Lal. Causal categories: relativistically interacting processes, preprint, (2011).

\bibitem{Com} M. Comeau. {\it Premonoidal $*$-Categories and Algebraic Quantum Field Theory}, 
Thesis, University of Ottawa, (2012).

\bibitem{DR1} S. Doplicher, J.E. Roberts, A new duality theory for compact groups,
Invent. math. 98, 157-218 (1989).

\bibitem{GLR} P. Ghez, R. Lima, and J. E. Roberts, $W^*$-categories.  Pacific J. Math. 120, pp. 
79-109, (1985).

\bibitem{Hal} H. Halvorson, Algebraic Quantum Field Theory, {\it Philosophy of
Physics}, edited by Jeremy Butterfield and John Earman, pp. 731-922, North-Holland (2006).

\bibitem{JS} V. Jones, V. Sunder. {\it Introduction to Subfactors}, Cambridge
University Press, (1997).

\bibitem{KR} R. Kadison, J. Ringrose, {\it Fundamentals of the Theory of Operator Algebras}, 
Graduate Studies in Mathematics, American Mathematical Society, (1997). 

\bibitem{CWM} S. Mac Lane, \emph{Categories for the Working Mathematician}, Second
Edition, Springer-Verlag, (2000). 

\bibitem{Pen} R. Penrose. {\em Techniques of Differential Topology in Relativity}.
CBMS-NSF Regional Conference Series in Applied Mathematics, (1972).

\bibitem{PR} J. Power, E. Robinson. Premonoidal categories and models of computation.
Mathematical Structures in Computer Science  7 pp. 453-468, (1997). 

\bibitem{Rob} J. Roberts. Lectures on Algebraic Quantum Field Theory. In: {\it The
Algebraic Theory of Superselection Sectors.} , ed. D. Kastler, pp. 1-112.
World ScientiÞc, (1990).

\bibitem{Sel} P. Selinger. Dagger compact closed categories and completely positive maps.
Electronic Notes in Theoretical Computer Science 170, pp. 139-163, (2007).
 
\bibitem{Sel2} P. Selinger. Finite-dimensional Hilbert spaces are complete for dagger compact closed categories.
{\sf preprint}, (2012).
 
\bibitem{Sun} V. Sunder, {\it An invitation to Von Neumann Algebras},
Springer-Verlag, (1987).

\bibitem{Will} D. Williams. {\it Crossed Products of $C^*$-algebras}, American
Mathematical Society, (2007).

\end{thebibliography}
\end{document}